\documentclass[11pt]{article}
\PassOptionsToPackage{obeyspaces}{url}
\usepackage[hidelinks]{hyperref}
\usepackage{bookmark}
\bookmarksetup{open,numbered,depth=5}
\usepackage{amsfonts}
\usepackage{amsmath}
\usepackage{amssymb, amsthm, enumitem}
\usepackage{breakurl}
\usepackage{tikz}
\usetikzlibrary{calc}

\usepackage{etoolbox} 
\patchcmd{\thebibliography}{\leftmargin\labelwidth}{\leftmargin\labelwidth\addtolength\itemsep{-0.1\baselineskip}}{}{}

\usepackage{mathtools}
\usepackage{xstring}

\author{Zichao Dong\thanks{Department of Mathematical Sciences, Carnegie Mellon University, Pittsburgh, PA 15213, USA\@. Supported in part by U.S.\ taxpayers through NSF CAREER grant DMS-1555149.} \and Zhuo Wu\thanks{Beijing International Center for Mathematical Research, Peking University, Beijing 100871, China. }}

\title{On the stability of graph independence number}
\date{}

\usepackage[nameinlink]{cleveref}

\newtheorem{theorem}{Theorem}
\newtheorem{lemma}[theorem]{Lemma}
\newtheorem{corollary}[theorem]{Corollary}

\newcommand*{\eqdef}{\stackrel{\mbox{\normalfont\tiny def}}{=}} 
\newcommand*{\veps}{\varepsilon}                                

\newcommand*{\N}{\mathbb{N}}                                    
\newcommand*{\Z}{\mathbb{Z}}                                    
\newcommand*{\cK}{\mathcal{K}}
\newcommand*{\cM}{\mathcal{M}}
\newcommand*{\cP}{\mathcal{P}}
\newcommand*{\cG}{\mathcal{G}}

\crefname{enumi}{step}{steps}
\crefname{part}{part}{parts}

\begin{document}
\maketitle

\begin{abstract}
	Let $G$ be a graph on $n$ vertices of independence number $\alpha(G)$ such that every induced subgraph of $G$ on $n-k$ vertices has an independent set of size at least $\alpha(G) - \ell$. What is the largest possible $\alpha(G)$ in terms of $n$ for fixed $k$ and $\ell$? We show that $\alpha(G) \le n/2 + C_{k, \ell}$, which is sharp for $k-\ell \le 2$. We also use this result to determine new values of the Erd\H{o}s--Rogers function. 
\end{abstract}

\section{Introduction}

\paragraph{Background.} All graphs considered here are finite, undirected, and simple. For graph $G$ and property $\cP$, the resilience of $\cP$ measures how much one should change $G$ in order to destroy $\cP$.  Assume $G$ has $\cP$, then the \emph{global resilience} of $\cP$ refers to the minimum number $r$ such that by removing $r$ edges from $G$ one can obtain a graph not having $\cP$, and the \emph{local resilience} refers to the minimum number $r$ such that by removing at each vertex at most $r$ edges one can obtain a graph not having $\cP$. For example, Tur\'{a}n's theorem (see \cite{turan}) characterizes the global resilience of having a $k$-clique in complete graphs, and Dirac's theorem (see \cite{dirac}) characterizes the local resilience of having a Hamiltonian cycle in complete graphs. Moreover, the local resilience of various properties is extensively studied in \cite{sudakov_vu}. 

Note that both global resilience and local resilience focus on removing edges. What about removing vertices? As far as we are aware of, this vertex-removal version of resilience is never discussed. To distinguish from removing edges, we shall always use the word stability when discussing removing vertices throughout this paper. 

To be more specific, we are going to study the resilience of graph independence number with respect to removing vertices. For vertices $v_1, \ldots, v_k$, denote the induced subgraph of $G$ on $V(G) \setminus \{v_1, \ldots, v_k\}$ by $G \setminus \{v_1, \ldots, v_k\}$. For non-negative integers $k > \ell \ge 0$, a graph $G(V, E)$ is called \emph{$(k, \ell)$-stable}, if for every $k$ vertices $v_1, \ldots, v_k$ of $G$, 
\[
\alpha(G \setminus \{v_1, \ldots, v_k\}) \ge \alpha(G) - \ell. 
\]
That is, the independence number $\alpha(G)$ drops by at most $\ell$ after removing any $k$ vertices from $V(G)$. 

This is related to a classical problem of Erd\H{o}s and Rogers. Extending the problem of studying Ramsey numbers, in \cite{erdos_rogers} they defined the function
\[
f_{s, s+t}(n) \eqdef \min \big\{ \max\{ |S|: \text{$S \subseteq V(G)$ and the induced subgraph $G[S]$ contains no $K_s$} \} \big\}, 
\]
where the minimum is taken over all $K_{s+t}$-free graphs $G$ on $n$ vertices. They established a lower bound on Ramsey number $R(k, \ell)$ by arguing that $f_{s, s+1}(n) \le n^{1-\veps(s)}$ for some positive constant $\veps(s)$. Throughout the years, the upper bounds and lower bounds on $f_{s, s+t}(n)$ for different pairs $(s, t)$ have been extensively studied (see \cite{bollobas_hind,krivelevich94,krivelevich95,sudakov05a,sudakov05b,dudek_rodl,wolfovitz,dudek_retter_rodl,gowers_janzer}). The case $t=s+1$ has received the most attention. The best known bounds (see \cite{dudek_mubayi} and \cite{dudek_retter_rodl}) in this case are
\[
\Omega \bigl( n^{1/2} (\log n)^{1/2} (\log \log n)^{-1/2}\bigr) \le f_{s, s+1}(n) \le O \bigl( n^{1/2} (\log n)^{4s^2} \bigr). 
\]
For a detailed survey on Erd\H{o}s--Rogers function, see \cite{dudek_rodl_survey}. For a recent study of an alternative local-global perspective on the Erd\H{o}s--Rogers function, see \cite{bucic_sudakov}.  

The study of Erd\H{o}s--Rogers function has focused on the case when $s, t$ are fixed and $n$ tends to infinity. Our results imply the exact value of $f_{s, s+t}(n)$ when $s > n/2$.

\paragraph{Results.} Unlike the previous work on the Erd\H{o}s--Rogers problem, in this paper we study the behavior of independence number on very large induced subgraphs of a graph. That is, given a $(k, \ell)$-stable graph $G(V, E)$, what can be said about $\alpha(G)$ in terms of $n \eqdef |V|$ as $n \to \infty$? Our main result is the following upper bound. 

\begin{theorem} \label{thm:pkl}
	Suppose $G(V, E)$ is a $(k, \ell)$-stable graph with $|V(G)| = n$ and $k > \ell \ge 0$, then 
	\begin{equation} \label{eqn:pkl}
		\alpha(G) \le \left\lfloor \frac{n-k+1}{2} \right\rfloor + \ell. 
	\end{equation}
\end{theorem}

This upper bound on $\alpha(G)$ determines a new class of Erd\H{o}s--Rogers function values. 

\begin{theorem} \label{thm:ERfunc}
	$f_{s, s+t}(n) = n-t$ for every integers $s, t, n$ such that $s > \lfloor \frac{n-t+1}{2} \rfloor$ and $s+t \le n+1$. 
\end{theorem}

To prove \Cref{thm:pkl}, the crucial step is to prove its special case when $(k, \ell) = (1, 0)$. This special case also has an interesting corollary. 

\begin{corollary} \label{cor:p1'}
	Suppose $G(V, E)$ is an $n$-vertex graph. If $\alpha(G \setminus \{v\}) = \alpha(G)$ holds for $m$ choices $v \in V$, then $\alpha(G) \le \lfloor n-\frac{m}{2} \rfloor$. 
\end{corollary}

Next, we explore the tightness of \Cref{thm:pkl}. We call a graph attaining equality in \Cref{thm:pkl} \emph{tight $(k, \ell)$-stable}. For example, every balanced complete bipartite graph is a tight $(1, 0)$-stable graph. 

Obviously, the complete graph $K_{k+1}$ is tight $(k, 0)$-stable for every $k$. However, it is harder to find tight $(k, \ell)$-stable graphs on larger vertex sets. The next theorems are in that direction. 

\begin{theorem} \leavevmode\label{thm:pkltight}
	Suppose $n > k > \ell \ge 0$. 
	\begin{enumerate}[label=(\roman*), ref=(\roman*)]
		\item \label[part]{pkltight:l+1}
		For every $\ell$, if $k=\ell+1$, then there exists an $n$-vertex tight $(k, \ell)$-stable graph for every $n$. 
		\item \label[part]{pkltight:l+2}
		For every $\ell$, if $k=\ell+2$, then there exists an $n$-vertex tight $(k, \ell)$-stable graph for every $n$. 
	\end{enumerate}
\end{theorem}

\begin{theorem} \label{thm:p3a}
	For every $\ell \ge 0$, if $k = \ell+3$, then there exists a sequence of $(\ell+3, \ell)$-stable graphs $G_m(V_m, E_m)$ with $|V_m| \to \infty$ such that 
	\[
	\alpha(G_m) = \frac{|V_m|}{2} - O\left(\sqrt{|V_m|}\right). 
	\]
\end{theorem}

\paragraph{Paper organization.} This paper is split into two parts. 

In the first part, we prove the upper bounds (\Cref{thm:pkl}). We then apply them to prove \Cref{thm:ERfunc} and \Cref{cor:p1'}. These occupy \Cref{sec:p1proof} through \Cref{sec:ERproof}. 

In the second part, we study the tightness of the upper bounds. We prove \Cref{thm:pkltight,thm:p3a} in \Cref{sec:pkltight,sec:p3aproof}. 

Finally, we devote \Cref{sec:t12} to a partial characterization of tight $(1, 0)$-stable graphs and tight $(2, 0)$-stable graphs. In particular, we prove tight lower and upper bounds on number of edges in tight $(1, 0)$-stable graphs.

\section{Proof of \texorpdfstring{\Cref{thm:pkl}}{Theorems~\ref{thm:p1}} assuming \boldmath $(k, \ell) = (1, 0)$} \label{sec:p1proof}

Let $G(V, E)$ be a $(1, 0)$-stable graph. Take a maximum-sized independent set $Y$ of $V$. Since $|Y| = \alpha(G)$, it suffices to show that Hall's condition (see, e.g., \cite{west}) holds from $Y$ to $V \setminus Y$. To be specific, for a set $A \subseteq V$, define 
\[
N(A) \eqdef \{u \in V \setminus A: \text{$u$ is a neighbor of some $v \in A$}\}. 
\]
Our goal is to prove that $|N(Y')| \ge |Y'|$ for every subset $Y' \subseteq Y$. 
	
Assume, for contradiction's sake, that $|N(Z)| < |Z|$ for some $Z \subseteq Y$. We may assume that $Z$ is minimal. Choose an arbitrary $z_0 \in Z$, which exists since $Z \neq \varnothing$. 
	
Because $\alpha(G) = \alpha(G \setminus \{z_0\})$, we can find another maximum-sized independent set $X$ with $z_0 \notin X$. Define $Z_1 = X \cap Z$ and $Z_2 = Z \setminus Z_1$. Put 
\[
U \eqdef \big( X \setminus N(Z_2) \big) \cup Z_2. 
\] 
We claim that $U$ is independent and that $|U| > \alpha(G)$, which would be a contradiction. 
	
First, we show that $U$ is independent. Note that both $X$ and $Z$ are independent, hence both $X \setminus N(Z_2)$ and $Z_2$ are independent. Since there is no edge between $Z_2$ and $X \setminus N(Z_2)$, we see that $U$ is independent. 
	
Next, we show that $U$ can also be written as $\big( X \setminus ( N(Z) \setminus N(Z_1) ) \big) \cup Z_2$. It suffices to show that $X \setminus ( N(Z) \setminus N(Z_1) ) = X \setminus N(Z_2)$. Suppose $x \in X \setminus (N(Z) \setminus N(Z_1))$ and $z \in Z_2$ are adjacent. Then $x \in N(Z_2)$ and hence $x \in N(Z)$. Since $Z_1 \subseteq X$ and $X$ is independent, there is no edge between $Z_1$ and $X$, hence $x \notin N(Z_1)$. Thus, $x \in X$ while $x \in N(Z) \setminus N(Z_1)$, which is a contradiction. We conclude that $X \setminus ( N(Z) \setminus N(Z_1) ) \subseteq X \setminus N(Z_2)$. Since $N(Z) \setminus N(Z_1) \subseteq N(Z_2)$, the opposite inclusion holds as well. Hence $U = \big( X \setminus ( N(Z) \setminus N(Z_1) ) \big) \cup Z_2$. 
	
Finally, we show that $|U| > |X| = \alpha(G)$. Note that $Z_2 \cap X = \varnothing$ by definition, so we have 
\[
|U| = |X \setminus (N(Z) \setminus N(Z_1))| + |Z_2| \ge |X| - |N(Z) \setminus N(Z_1)| + |Z_2|. 
\]
It suffices to show that $|N(Z) \setminus N(Z_1)| < |Z_2|$. Since $z_0 \in Z \setminus X$, so $Z_1 \subsetneq Z$, hence $|N(Z_1)| \ge |Z_1|$ by the minimality assumption of $Z$. Thus, 
\[
|N(Z) \setminus N(Z_1)| = |N(Z)| - |N(Z_1)| \le |N(Z)| - |Z_1| < |Z| - |Z_1| = |Z_2|, 
\]
hence $|U| > \alpha(G)$. 
	
This is the promised contradiction. Hence the proof is complete. \qed

\section{Proofs of \texorpdfstring{\Cref{thm:pkl}}{Theorems~\ref{thm:pkl}} and \texorpdfstring{\Cref{cor:p1'}}{Theorems~\ref{cor:p1'}}} \label{sec:pklproof}

In this section, we derive \Cref{thm:pkl} and \Cref{cor:p1'} from what we proved in \Cref{sec:p1proof}. 

\begin{proof}[Proof of \Cref{thm:pkl}.]
	The proof is by induction on $\ell$ and $k$. 
			
	Suppose $\ell = 0$, and we are going to show that 
	\begin{equation} \label{eqn:pk}
	\alpha(G) \le \left\lfloor \frac{n-k+1}{2} \right\rfloor. 
	\end{equation}
	\eqref{eqn:pk} holds when $k=1$ since we proved it in \Cref{sec:p1proof}. Assume \eqref{eqn:pk} is established for $k-1$, and suppose $G(V, E)$ is a $(k, 0)$-stable graph with $|V|=n$. For any vertex $v \in V$, by definition we have that $G \setminus \{v\}$ is $(k-1, 0)$-stable with $\alpha(G \setminus \{v\}) = \alpha(G)$, so the induction hypothesis tells us that 
	\[
	\alpha(G) = \alpha(G \setminus \{v\}) \le \left\lfloor \frac{(n-1)-(k-1)+1}{2} \right\rfloor = \left\lfloor \frac{n-k+1}{2} \right\rfloor, 
	\]
	which concludes the inductive proof of \eqref{eqn:pk}. 
			
	Assume \eqref{eqn:pkl} is established for $\ell-1$, and suppose $G(V, E)$ is a $(k, \ell)$-stable graph with $|V|=n$. If there exists $v_1, \ldots, v_{k-1} \in V$ such that $\alpha(G') = \alpha(G)-\ell$, here $G' \eqdef G \setminus \{v_1, \ldots, v_{k-1}\}$. Then by definition $G'$ is a $(1, 0)$-stable graph, hence what we proved in \Cref{sec:p1proof} tells us that 
	\[
	\alpha(G) \le \alpha(G') + \ell = \left\lfloor \frac{n-(k-1)}{2} \right\rfloor + \ell = \left\lfloor \frac{n-k+1}{2} \right\rfloor + \ell.
	\]
	Otherwise, $\alpha(G)$ drops by at most $\ell-1$ after removing any $(k-1)$-element subset of $V$, and by definition $G$ is $(k-1, \ell-1)$-stable, hence the induction hypothesis tells us that 
	\[
	\alpha(G) \le \left\lfloor \frac{n-(k-1)+1}{2} \right\rfloor + (\ell-1) = \left\lfloor \frac{n-k}{2} \right\rfloor + \ell. 
	\]
	Putting these together gives us that \eqref{eqn:pkl} holds for $k+1$, which concludes the inductive proof. 
\end{proof}

\begin{proof}[Proof of \Cref{cor:p1'}]
	Let $S$ be the set of vertices $v$ satisfying $\alpha(G \setminus \{v\}) = \alpha(G)$. By the assumption $|S| \ge m$. Let $\overline{S} = V(G) \setminus S$ be the complement of $S$. A key observation is that every maximum-sized independent set of $G$ contains $\overline{S}$. 
	
	Define 
	\[
	S_1 \eqdef \{v \in S: N(v) \cap \overline{S} = \varnothing\}, \qquad S_2 \eqdef \{v \in S: N(v) \cap \overline{S} \neq \varnothing\}. 
	\]
	Let $G_1$ be the induced subgraph of $G$ on $S_1$, and we claim that $\alpha(G_1 \setminus \{v\}) = \alpha(G_1)$ for every vertex $v$ of $S_1$. Indeed, for $v \in S_1$, let $X$ be an independent set of size $\alpha(G)$ in $G \setminus \{v\}$. Then $X_1 \eqdef X \cap S_1$ has to be an independent set of $G_1$ of size $\alpha(G_1)$, for otherwise there exists some independent set $X_1'$ of $G_1$ such that $|X_1'| > |X_1|$, hence $X' \eqdef X_1' \cup \overline{S}$ is an independent set in $G$ with $|X'| > \alpha(G)$, which is impossible. Thus, the existence of $X_1$ verifies our claim. 
	
	Since every maximum-sized independent set of $G$ contains $\overline{S}$, we have $\alpha(G) = |\overline{S}| + \alpha(G_1)$, hence it follows from what we proved in \Cref{sec:p1proof} that 
	\[
	\alpha(G) = |\overline{S}| + \alpha(G_1) \le n - |S| + |S_1|/2 \le n - |S| + |S|/2 \le n - m/2. \qedhere
	\]
\end{proof}

We remark that the bound in \Cref{cor:p1'} is tight, as witnessed by a disjoint union of a balanced complete bipartite graph $K_{\lceil m/2 \rceil, \lceil m/2 \rceil}$ and an independent set of size $n-2\lceil m/2 \rceil$ (possibly empty) when $m<n$. For the case $m=n$, see \Cref{sec:pkltight}. 

\section{Proof of \texorpdfstring{\Cref{thm:ERfunc}}{Theorems~\ref{thm:ERfunc}}} \label{sec:ERproof}

By considering the graph complements, there is no difference between cliques and independent sets in the definition of $f_{s, s+t}$. For convenience we replace all cliques by independent sets. 

Define $\cG_n^m \eqdef \{G: \text{$G$ is a graph on $n$ vertices with $\alpha(G) \le m-1$}\}$. Our goal is to verify that 
\begin{itemize}
	\item on one hand, there is $H \in \cG_n^{s+t}$ such that for every $t-1$ vertices $u_1, \ldots, u_{t-1}$ of $H$, 
	\[
	\alpha(H \setminus \{u_1, \ldots, u_{t-1}\}) \ge s; 
	\]
	\item on the other hand, for every $G \in \cG_n^{s+t}$, there are $t$ vertices $v_1, \ldots, v_t$ of $G$ such that 
	\[
	\alpha(G \setminus \{v_1, \ldots, v_t\}) \le s-1. 
	\]
\end{itemize}

The first part is seen by considering any $n$-vertex graph $H$ with $\alpha(H) = s+t-1$. Such an $H$ exists because $s+t \le n+1$. Note that the second part is trivial when $\alpha(G) \le s-1$, so we assume that $\alpha(G) = s-1+m$ for some positive integer $m \le t$. The key is to observe that $G$ is not $(t, m-1)$-stable. Indeed, if $G$ were $(t, m-1)$-stable, then \Cref{thm:pkl} would imply that 
\[
s-1+m = \alpha(G) \le \left\lfloor \frac{n-t+1}{2} \right\rfloor + (m-1), 
\]
hence $s \le \lfloor \frac{n-t+1}{2} \rfloor$, which is a contradiction. Now that $G$ is not $(t, m-1)$-stable, there are $t$ vertices $v_1, \ldots, v_t$ of $G$ such that 
\[
\alpha(G \setminus \{v_1, \ldots, v_t\}) \le \alpha(G) - m = s-1, 
\]
which concludes the proof of the second part. The proof of \Cref{thm:ERfunc} is complete. \qed

\section{Proof of \texorpdfstring{\Cref{thm:pkltight}}{Theorems~\ref{thm:pkltight}}} \label{sec:pkltight}

With the help of the lemma below, the general $(k, \ell)$-stability in \Cref{thm:pkltight} can be reduced to the simpler $(k-\ell, 0)$-stability. 

\begin{lemma} \label{lem:pkllemma}
	Suppose $n>k>\ell\ge0$. If $G$ is an $n$-vertex tight $(k, \ell)$-stable graph, then $G \sqcup K_1$, the disjoint union of $G$ and an isolated vertex, is an $(n+1)$-vertex tight $(k+1, \ell+1)$-stable graph. 
\end{lemma}

\begin{proof}
	Note that $G$ is tight $(k, \ell)$-stable. Hence
	\[
	\alpha(G \sqcup K_1) = \alpha(G)+1 = \left(\left\lfloor \frac{n-k+1}{2} \right\rfloor + \ell\right) + 1 = \left\lfloor \frac{(n+1)-(k+1)+1}{2} \right\rfloor + (\ell+1), 
	\]
	which attains equality in \Cref{thm:pkl}. So it suffices to prove that $G \sqcup K_1$ is $(k+1, \ell+1)$-stable. 
	
	Let $v'$ be the isolated vertex in $G \sqcup K_1$, and suppose a subset of $k+1$ vertices $\{v_1, \ldots, v_{k+1}\}$ are removed from $V(G \sqcup K_1)$. If $v'$ is removed, we may assume $v_{k+1}=v'$, hence by the stability of $G$, 
	\[
	\alpha(G \sqcup K_1 \setminus \{v_1, \ldots, v_{k+1}\}) = \alpha(G \setminus \{v_1, \ldots, v_k\}) \ge \alpha(G) - \ell = \alpha(G \sqcup K_1) - (\ell+1). 
	\]
	Otherwise, $v_1, \ldots, v_{k+1} \in V(G)$, then by the stability of $G$, 
	\begin{align*}
	\alpha(G \sqcup K_1 \setminus \{v_1, \ldots, v_{k+1}\}) &= \alpha \big( (G \setminus \{v_1, \ldots, v_k\}) \setminus \{v_{k+1}\} \big) + 1 \\
	&\ge \alpha(G \setminus \{v_1, \ldots, v_k\}) \ge \alpha(G) - \ell = \alpha(G \sqcup K_1) - (\ell+1). 
	\end{align*}
	Thus, $G \sqcup K_1$ is $(k+1, \ell+1)$-stable, and the proof is complete. 
\end{proof}

Evidently, \Cref{lem:pkllemma} helps reduce the proof of \Cref{thm:pkltight}\ref{pkltight:l+1}, \Cref{thm:pkltight}\ref{pkltight:l+2}, and \Cref{thm:p3a} in general to the proof of their special cases when $(k, \ell) = (1, 0)$, $(2, 0)$, and $(3, 0)$. 

\begin{proof}[Proof of \Cref{thm:pkltight}]
	For $n \ge 2$, denote by $\cK_n$ the complete bipartite graph $K_{n/2, n/2}$ when $n$ is even, or the complete tripartite graph $K_{\lfloor n/2 \rfloor, \lfloor n/2 \rfloor, 1}$ when $n$ is odd. Then $\alpha(\cK_n) = \lfloor n/2 \rfloor$. After removing any one vertex $v$, at least one $\lfloor n/2 \rfloor$-vertex part is left untouched, so $\alpha(\cK_n \setminus \{v\}) = \lfloor n/2 \rfloor = \alpha(\cK_n)$. Hence $\cK_n$ is always $(1, 0)$-stable. By applying \Cref{lem:pkllemma}, we see that \Cref{thm:pkltight}\ref{pkltight:l+1} holds. 
	
	As for \Cref{thm:pkltight}\ref{pkltight:l+2}, it suffices to construct an $n$-vertex $(2, 0)$-stable graph for every $n \ge 3$. 
	
	Consider \emph{cycles} $C_n$ and \emph{wheels} $W_n$. Here an \emph{$n$-wheel} refers to an $(n-1)$-cycle $C_{n-1}$ with another vertex connected to every vertex of the cycle. For example, $C_5$ and $W_6$ are shown in the figure below: 
	\[ 
	\begin{tikzpicture} 
		\fill (-3, 0) circle (2pt); 
		\fill (-2.5, -1) circle (2pt); 
		\fill (-1, 0) circle (2pt); 
		\fill (-1.5, -1) circle (2pt); 
		\fill (-2, 0.5) circle (2pt); 
		\draw (-3, 0) -- (-2.5, -1); 
		\draw (-2.5, -1) -- (-1.5, -1); 
		\draw (-1.5, -1) -- (-1, 0); 
		\draw (-1, 0) -- (-2, 0.5); 
		\draw (-2, 0.5) -- (-3, 0); 
		\fill (1, 0) circle (2pt); 
		\fill (1.5, -1) circle (2pt); 
		\fill (3, 0) circle (2pt); 
		\fill (2.5, -1) circle (2pt); 
		\fill (2, 0.5) circle (2pt); 
		\fill (2, -0.3) circle (2pt); 
		\draw (1, 0) -- (1.5, -1); 
		\draw (1.5, -1) -- (2.5, -1); 
		\draw (2.5, -1) -- (3, 0); 
		\draw (3, 0) -- (2, 0.5); 
		\draw (2, 0.5) -- (1, 0); 
		\draw (2, -0.3) -- (1, 0); 
		\draw (2, -0.3) -- (1.5, -1); 
		\draw (2, -0.3) -- (2.5, -1); 
		\draw (2, -0.3) -- (3, 0); 
		\draw (2, -0.3) -- (2, 0.5); 
		\node at (0, -1.5) {\textbf{Figure 1: }Graphs $C_5$ and $W_6$.};
	\end{tikzpicture} 
	\]
	We claim that $C_n$ is tight $(2, 0)$-stable when $n$ is odd, and $W_n$ is tight $(2, 0)$-stable when $n$ is even. 
	
	Denote the $n$-vertex path graph as $P_n$, evidently $\alpha(P_n) = \lceil n/2 \rceil = \lfloor (n+1)/2 \rfloor$. 
	When $n$ is odd, obviously $\alpha(C_n) = (n-1)/2$. Suppose two disjoint paths $P_a$ and $P_b$ are left after vertices $u, v$ being removed from $C_n$ (note that $a, b$ can be zero), then $a+b = n-2$, and 
	\[
	\alpha(C_n \setminus \{u, v\}) = \left\lfloor \frac{a+1}{2} \right\rfloor + \left\lfloor \frac{b+1}{2} \right\rfloor = \frac{a+b+1}{2} = \frac{n-1}{2} = \alpha(C_n), 
	\]
	since $a$ and $b$ are of different parity. So $C_n$ is tight $(2, 0)$-stable. When $n$ is even, obviously $\alpha(W_n) = \alpha(C_{n-1})$. Suppose two vertices $u, v$ are removed from $W_n$, then at most two vertices are removed from the induced subgraph $C_{n-1}$. Since $C_{n-1}$ is $(2, 0)$-stable, we see that 
	\[
	\alpha(C_{n-1}) = \alpha(C_{n-1} \setminus \{u, v\}) \le \alpha(W_n \setminus \{u, v\}) \le \alpha(W_n) = \alpha(C_{n-1}), 
	\]
	so $\alpha(W_n) = \alpha(C_{n-1}) = n/2-1$, and hence $W_n$ is tight $(2, 0)$-stable. 
	
	Thus, an $n$-vertex tight $(2, 0)$-stable graph always exists. By applying \Cref{lem:pkllemma}, we are done. 
\end{proof}

It is worth mentioning that no $6$-vertex tight $(3, 0)$-stable graph exists. The proof is simple. Assume $G$ is a $6$-vertex tight $(3, 0)$-stable graph, then $\alpha(G) = \lfloor(6-2)/2\rfloor = 2$, so there is no triangle subgraph $K_3$ in the complement graph $\overline{G}$. Also, if three vertices of $G$ form a triangle, by removing the other $3$ vertices the independence number drops to $1$, which contradicts with the $(3, 0)$-stability. Thus, there is no triangle subgraph $K_3$ in $G$. Now we reach a contradiction with the fact that Ramsey number $R(3, 3) = 6$, so no $6$-vertex tight $(3, 0)$-stable graph exists.

\section{Proof of \texorpdfstring{\Cref{thm:p3a}}{Theorems~\ref{thm:p3a}}} \label{sec:p3aproof}

We are going to construct a sequence of $(3, 0)$-stable graphs $G_m(V_m, E_m)$ with $|V_m| \to \infty$ such that 
\[
\alpha(G_m) = |V_m|/2 - O \left(\sqrt{|V_m|}\right). 
\]
Note that the existence of such a sequence directly implies \Cref{thm:p3a} by \Cref{lem:pkllemma}. 

Suppose $m \in \N$ ($m \ge 3$ for technical reasons). Take 
\[
V_m \eqdef \Z/(2m^2+2m)\Z, \qquad E_m \eqdef \{(i, j): |i-j| \equiv m \text{ or } m+1 \pmod {2m^2+2m}\}, 
\]
we claim that $G_m \eqdef (V_m, E_m)$ is $(3, 0)$-stable with $\alpha(G_m) = m^2$. Since 
\[
m^2 = (2m^2+2m)/2 - O\left(\sqrt{2m^2+2m}\right), 
\]
\Cref{thm:p3a} is shown once we verify the claim. 

\paragraph{First, we prove that \boldmath $\alpha(G_m) = m^2$.} It is easily seen that the vertices 
\[
\begin{array}{ccccc}
	0, & 2m+1, & 4m+2, & \cdots, & (m-1)(2m+1), \\
	1, & 2m+2, & 4m+3, & \cdots, & (m-1)(2m+1)+1, \\
	2, & 2m+3, & 4m+4, & \cdots, & (m-1)(2m+1)+2, \\
	\vdots & \vdots & \vdots & \ddots & \vdots \\
	m-1, & 3m, & 5m+1, & \cdots, & (m-1)(2m+1)+m-1 
\end{array}
\]
form an independent set of size $m^2$ in $G_m$. Hence, it suffices to show that $\alpha(G_m) \le m^2$. 

We prove this by contradiction. Suppose $X \subset V_m$ is an independent set of size $m^2+1$. Partition the vertices into $m$ groups of size $2m+2$ each as follows: 
\begin{align*}
	V_0 &\eqdef \{0, m+0, 2m+0, \cdots, (2m+1)m+0\}, \\
	V_1 &\eqdef \{1, m+1, 2m+1, \cdots, (2m+1)m+1\}, \\
	V_2 &\eqdef \{2, m+2, 2m+2, \cdots, (2m+1)m+2\}, \\
	&\ \ \vdots \\
	V_{m-1} &\eqdef \{m-1, 2m-1, 3m-1, \cdots, (2m+2)m-1\}. 
\end{align*}
Note that $|X \cap V_j| \le m+1$ since every two consecutive points (including the first and the last) are adjacent in $G_m$. Thus, by pigeonhole principle $|X \cap V_j| = m+1$ for some $j$. Due to the obvious translation symmetry in $G_m$, we may assume without loss of generality that 
\[X \cap V_0 = \{0, 2m, 4m, \ldots, 2m^2\}. \]
Since $X$ is an independent set, from the construction of $E_m$ we see that 
\[
\begin{array}{ccccccc}
	m, &3m, &5m, &\ldots, &(2m+1)m &\notin &X, \\
	m+1, &3m+1, &5m+1, &\ldots, &(2m+1)m+1 &\notin &X, \\
	m-1, &3m-1, &5m-1, &\ldots, &(2m+1)m-1 &\notin &X.
\end{array}
\]
The key is to pair up every other vertex of $G_m$ as follows: 
\[
\begin{array}{cccc} 
	\{1, m+2\}, & \{2m+1, 3m+2\}, & \cdots, & \{(2m \cdot m + 1, (2m+1)m + 2)\}, \\
	\{2, m+3\}, & \{2m+2, 3m+3\}, & \cdots, & \{(2m \cdot m + 2, (2m+1)m + 3)\}, \\
	\vdots & \vdots & \ddots & \vdots \\
	\{m-2, 2m-1\}, & \{3m-2, 4m-1\}, & \cdots, & \{(2m+1)m - 2, (2m+2)m - 1\}.
\end{array}
\]
Note that the vertices in each pair are adjacent in $G_m$. Hence they are not in $X$ simultaneously. So 
\[
|X| \le |X \cap V_0| + \text{\#pairs} = (m+1) + (m-2)(m+1) = m^2-1, 
\]
which contradicts with $|X| = m^2+1$. Thus, $\alpha(G) = m^2$, as desired. 

\paragraph{Then, we show that \boldmath $G_m$ is $(3, 0)$-stable.} 
Think about $V(G_m)$ as $2m^2 + 2m$ evenly distributed points on a circle. Here two points are adjacent in $G_m$ if and only if they form an interval of length $m$ or $m+1$ on the circle. We are going to prove that no matter which $3$ points are deleted, we can still pick $m^2$ points such that no two of them are adjacent to each other in $G_m$. 

Suppose $3$ points $u, v, w$ (clockwise in this order) are deleted, and the circle is cut into three pieces with $a, b, c$ consecutive points in $u \frown v, v \frown w, w \frown u$, respectively ($a, b, c$ can be zero), then $a+b+c = 2m^2+2m-3$. Here $u \frown v$ refers to the clockwise arc from $u$ to $v$ on the circle. The idea is to chop these pieces into many \emph{normal groups}, $(2m+1)$-consecutive-point groups, and three small consecutive-point groups, then to pick some $m$ consecutive points out of each normal group, and to deal with the remaining small groups carefully. Set 
\[
a' \eqdef a \bmod {2m+1}, \qquad b' \eqdef b \bmod {2m+1}, \qquad c' \eqdef c \bmod {2m+1}. 
\]

\underline{Case 1. $a'+b'+c' = m-3$.} Pick the first $m$ points clockwise in each normal group. We have picked $m^2$ points in total, and every two of them are either among some $m$ consecutive points or are at least $m+2$ apart from each other, so these points form an independent set of size $m^2$. 

\underline{Case 2. $a'+b'+c' = 3m-2$.} By symmetry we may assume without loss of generality that $a' \ge m$. Clockwise chop $u \frown v$ into $a', 2m+1, \ldots, 2m+1$-point groups, chop $v \frown w$ into $2m+1, \ldots, 2m+1, b'$-point groups, chop $w \frown u$ into $2m+1, \ldots, 2m+1, c'$-point groups. Then we pick points as below: 
\begin{itemize}
	\item Pick the first $m$ points clockwise in the $a'$-point group (recall that $a' \ge m$). 
	\item Starting from the last picked point clockwise in the $a'$-point group and moving clockwise on the circle, we pick the first \emph{available} $m$ points from each normal group (\emph{available} means forming an interval of length at least $m+2$ on the circle from the last point that has been picked). 
\end{itemize}
So far we have picked $m^2$ points. By the way we picked these points, 
\begin{itemize}
	\item the consecutive $a'-m$ points clockwise before $v$ are unpicked; 
	\item the consecutive $b'+(a'-m)+1$ points clockwise before $w$ are unpicked; 
	\item the consecutive $c'+\min\{b'+(a'-m)+2, m+1\}$ points clockwise before $u$ are unpicked. 
\end{itemize}
According to the way we picked these points, the open interval between consecutive picked points including $v$ and the open interval between consecutive picked points including $w$ (possibly the same) both contain at least $m+1$ points. As for the open interval between consecutive picked points including $u$, it also contains at least $m+1$ points, as $a'+b'+c' = 3m-2$ implies that
\[
c'+\min\{b'+(a'-m)+2, m+1\}+1 \ge m+1. 
\]
Thus, the $m^2$ points we picked form an independent set. 

\underline{Case 3. $a'+b'+c' = 5m-1$.} The fact $a', b', c' \le 2m$ implies that $m-1 \le a', b', c' \le 2m$. By symmetry we may assume $a' = \min\{a', b', c'\}$. Clockwise chop $u \frown v$ into $2m+1, \ldots, 2m+1, a'$-point groups, chop $v \frown w$ into $b', 2m+1, \ldots, 2m+1$-point groups, chop $w \frown u$ into $2m+1, \ldots, 2m+1, c'$-point groups. Then we pick points as below: 
\begin{itemize}
	\item Starting from $u$ and moving clockwise on $u \frown v$, we pick the $(m+1)$'th through the $2m$'th point from each normal group. 
	\item Starting from $v$ and moving clockwise, we pick the first $m$ points from the $b'$-point group (recall that $b' \ge m$). Starting from the last picked point and moving clockwise on $v \frown w$, we pick the first available $m$ points from each normal group. 
	\item Starting from $u$ and moving counterclockwise, we pick the first $m$ points from the $c'$-point group (recall that $c' \ge m$). Starting from the last picked point and moving counterclockwise on $w \frown u$, we pick the first available $m$ points from each normal group. 
\end{itemize}
So far we have picked $m^2$ points. By the way we picked these points, 
\begin{itemize}
	\item the consecutive $m$ points clockwise after $u$ are unpicked (because $a' \ge m-1$); 
	\item the consecutive $a'+1$ points clockwise before $v$ are unpicked; 
	\item the consecutive $b'-m$ points clockwise before $w$ are unpicked, and the consecutive $c'-m$ points clockwise after $w$ are unpicked. 
\end{itemize}
Since $a' \ge m-1$, the open interval between consecutive picked points including $u$ and the open interval between consecutive picked points including $v$ (possibly the same) both contain at least $m+1$ points. As for the open interval between consecutive picked points including $w$, it also contains at least $m+1$ points, as $a' = \min\{a', b', c'\}$ implies that $a' < 2m$, hence
\[
(b'-m)+(c'-m)+1 = (5m-1) - a' - 2m + 1 \ge m+1. 
\]
Thus, the $m^2$ points we picked form an independent set. 

By combining the three cases, we conclude that $G_m$ is $(3, 0)$-stable. 

The claim is seen by combining the above, and the proof is complete. \qed

\bigskip

We remark that \Cref{thm:pkl} is also asymptotically tight when $k = \ell+4$. In fact, every graph $G_m$ with $V_m \eqdef \Z/(2m^2+2m+1)\Z$ and $E_m \eqdef \{(i, j): |i-j| \equiv m \text{ or } m+1 \pmod {2m^2+2m+1}\}$ is~$(4, 0)$-stable, and the proof is basically the same as the proof of \Cref{thm:p3a} but slightly subtler. 

\section{More on tight \boldmath $(1, 0)$-stable and tight $(2, 0)$-stable graphs} \label{sec:t12}

As we have seen in the proof of \Cref{thm:pkltight}, there is an $n$-vertex $(1, 0)$-stable graph for every $n \ge 2$ and an $n$-vertex $(2, 0)$-stable graph for every $n \ge 3$. In this section, we care about the \emph{uniqueness} of such tight stable graphs. That is, for some given $n$, we want to know whether every $n$-vertex tight $(1, 0)$-stable graph is isomorphic to $\cK_n$ and whether every $n$-vertex $(2, 0)$-stable graph is isomorphic to $C_n$ or $W_n$. 

For tight $(1, 0)$-stable graphs, the uniqueness does not hold (assume $n \ge 4$). For example, denote by $\cM_n$ the disjoint union of $n/2$ edges when $n$ is even, and the disjoint union of $\lfloor n/2 \rfloor - 1$ edges and one triangle when $n$ is odd. Evidently, the graph $\cM_n$ is tight $(1, 0)$-stable for every $n$. Moreover, every $n$-vertex graph that is sandwiched between $\cM_n$ and $\cK_n$ is also tight $(1, 0)$-stable. 

Here we would like to mention that $n$-vertex ($n \ge 3$) path $P_n$ is another class of tight $(2, 1)$-stable graphs, which can be easily checked. This implies that there exist tight $(2,1)$-stable graphs other than those constructed in \Cref{lem:pkllemma}. 

For $(2, 0)$-stable graphs, the uniqueness does not hold when $n$ is even (assume $n \ge 6$). One can check that when 
\[
V_k \eqdef \Z/2k\Z, \qquad E_k \eqdef \{(i, j): |i-j| \equiv 1 \text{ or } k \pmod {2k}\}, 
\]
$G_k \eqdef (V_k, E_k)$ is tight $(2, 0)$-stable for every $k \ge 3$. However, we are unaware of any $n$-vertex tight $(2, 0)$-stable graph other than the odd cycle $C_n$ when $n$ is odd. We doubt that $C_n$ is the only $n$-vertex tight $(2, 0)$-stable graph. Note that $C_n$, $W_n$, and the tight $(2, 0)$-stable graphs constructed above all have Hamiltonian cycles inside. We also suspect that every tight $(2, 0)$-stable graph contains a Hamiltonian cycle. 

\paragraph{Characterization of tight \boldmath $(1, 0)$-stable graphs.} One might wonder whether a tight $(1, 0)$-stable graph $G$ (especially when $|V(G)|$ is even) is always sandwiched between a perfect matching and a complete bipartite graph (i.e. $\mathcal{M}_n \subseteq G \subseteq \mathcal{K}_n$). In fact, $\mathcal{M}_n \subseteq G$ is already shown by Hall's theorem in the proof in \Cref{sec:p1proof}, while $G \subseteq \mathcal{K}_n$ does not necessarily happen. An example is the graph below. 
\[
\begin{tikzpicture}
	\fill (-1, 1) circle (2pt);
	\fill (-2, 0) circle (2pt);
	\fill (-1, -1) circle (2pt);
	\fill (1, 1) circle (2pt);
	\fill (2, 0) circle (2pt);
	\fill (1, -1) circle (2pt);
	\draw (-1, 1) -- (1, 1);
	\draw (-2, 0) -- (2, 0);
	\draw (-1, -1) -- (1, -1);
	\draw (-1, 1) -- (-2, 0);
	\draw (-2, 0) -- (-1, -1);
	\draw (-1, -1) -- (-1, 1);
	\node at (0, -1.5) {\textbf{Figure 2: }A tight $(1,0)$-stable graph $G$ with $G \not\subseteq \cK_n$.};
\end{tikzpicture}
\]

Nevertheless, we can still prove tight upper and lower bounds on number of edges: 

\begin{theorem}
	Suppose $G(V, E)$ is a tight $(1, 0)$-stable graph with $|V(G)| = n$, then
	\[
	|E(\mathcal{M}_n)| \le |E(G)| \le |E(\mathcal{K}_n)|. 
	\]
\end{theorem}

\begin{proof} 
	We consider the cases when $n$ is even and when $n$ is odd separately. 
	
	\underline{Case 1. $n$ is even.} Suppose $n=2k$, then $\alpha(G) = k$ implies that the complement graph $\overline{G}$ contains no $(k+1)$-clique. Hence Tur\'{a}n's theorem implies that 
	\[
	|E(\overline{G})| \le |E( K_{2, \ldots, 2})| = 2k^2 - 2k. 
	\]
	Here $K_{2, \ldots, 2}$ refers to the complete $k$-partite graph in which each part contains $2$ vertices. So
	\[
	|E(G)| \ge \binom{2k}{2} - (2k^2-2k) = k = |E(\mathcal{M}_n)|, 
	\]
	which gives the desired lower bound. 
	
	For the upper bound, if there is a vertex $v$ with $\deg(v) \ge k+1$, then $v$ is not contained in any maximum-sized independent set since $\alpha(G) = k$. Hence $G \setminus \{v\}$ is $(1, 0)$-stable as well. However, $\alpha(G \setminus \{v\}) = k$, which contradicts \Cref{thm:pkl}. Thus, every vertex of $G$ is of degree at most $k$, hence 
	\[ |E(G)| \le \frac{1}{2} \cdot 2k \cdot k = k^2 = |E(\mathcal{K}_n)|, \]
	which gives the desired upper bound. 
	
	\underline{Case 2. $n$ is odd.} Suppose $n=2k+1$, then $\alpha(G) = k$ implies that the complement graph $\overline{G}$ contains no $(k+1)$-clique. Hence Tur\'{a}n's theorem implies that 
	\[
	|E(\overline{G})| \le |E( K_{2, \ldots, 2, 3})| = 2k^2 - 2. 
	\]
	Here $K_{2, \ldots, 2, 3}$ refers to the complete $k$-partite graph in which $k-1$ parts contain $2$ vertices in each and $1$ part contains $3$ vertices. So
	\[ |E(G)| \ge \binom{2k+1}{2} - (2k^2-2) = k+2 = |E(\mathcal{M}_n)|, \]
	which gives the desired lower bound. 
	
	For the upper bound, if there is a vertex $v$ with $\deg(v) \ge k+2$, then $v$ is not contained in any maximum-sized independent set since $\alpha(G) = k$. Hence $G \setminus \{v\}$ is $(1, 0)$-stable as well. So $|E(G \setminus \{v\})| \le k^2$ by what we just proved, hence
	\[ |E(G)| \le k^2 + \deg(v) \le k^2 + 2k = |E(\mathcal{K}_n)|. \]
	Otherwise, every vertex in $G$ is of degree at most $k+1$, hence 
	\[ |E(G)| \le \frac{1}{2} (2k+1) (k+1) \le k^2 + 2k = |E(\mathcal{K}_n)|. \]
	
	The proof is done by combining the two cases. 
\end{proof}

\section{Open problems} \label{sec:op}

\begin{enumerate}
	\item Is the odd cycle $C_n$ the only $n$-vertex tight $(2, 0)$-stable graph for odd $n \ge 3$? By exhaustive computer search, we verified this for $n = 3, 5, 7, 9, 11$. For even $n$, we suspect that every connected $n$-vertex tight $(2, 0)$-stable graph contains a Hamiltonian cycle. 
	
	\item \Cref{thm:pkl} is not always tight. For example, as we remarked in \Cref{sec:pkltight}, there exists no $6$-vertex tight $(3, 0)$-stable graph. It would be interesting to find the largest independence number of $n$-vertex $(k, \ell)$-stable graphs for every $n$, $k$, and $\ell$. 
	
	\item We do not know if \Cref{thm:pkl} is asymptotically tight, i.e., if there exist $n$-vertex $(k, \ell)$-stable graphs with independence number $n/2-o(n)$, for every $k$ and $\ell$. 
	
	\textbf{Remark.} This problem is answered positively by Theorem 1.4 in \cite{alon}. 
\end{enumerate}

\section*{Acknowledgments} 

We are grateful to Boris Bukh for proposing this problem to us, and for many helpful suggestions on writing. We thank Minghui Ouyang for help with computer programming, and for beneficial discussions. The second author is grateful to the support of Beijing International Center for Mathematical Research (BICMR) during his visit to BICMR in fall 2020 and spring 2021. We thank two anonymous referees for valuable feedback on the earlier versions of the paper.

\bibliographystyle{plain}
\bibliography{stable_indep}

\end{document}